\documentclass[12pt]{article}
\usepackage{amsmath, amsthm, amssymb}
\DeclareSymbolFont{AMSb}{U}{msb}{m}{n}
\DeclareMathSymbol{\R}{\mathbin}{AMSb}{"52}
\DeclareMathSymbol{\Z}{\mathbin}{AMSb}{"5A}

\begin{document}
\newtheorem{thm}{Theorem}[section]
\newtheorem{lemma}[thm]{Lemma}
\newtheorem{conj}[thm]{Conjecture}
\renewcommand{\qedsymbol}{}

\title{A Summary of Results and Problems Related to the Caccetta-H\"{a}ggkvist Conjecture}
\author{Blair D. Sullivan \\ Dept. of Mathematics, Princeton University \\ bdowling@princeton.edu}
\maketitle

\section{Introduction}
This paper is an attempt to survey the current state of our
knowledge on the Caccetta-H\"{a}ggkvist conjecture and related questions.
In January 2006 there was a workshop hosted by the American Institute of
Mathematics in Palo Alto, on the Caccetta-H\"{a}ggkvist conjecture, and this
paper partly originated there, as a summary of the open problems
and partial results presented at the workshop. Our thanks to the
many participants who helped with this paper.

\subsection*{Notations \& Definitions}
\begin{itemize}
\item A {\it graph} $G = (V(G), E(G))$ is a collection of vertices 
$V(G)$ and edges $E(G)$, where 
each edge is an unordered pair of distinct vertices ${u,v}$. 
\item A {\it digraph} $D = (V(D), E(D))$ is a collection of vertices 
$V(D)$ and edges $E(D)$, along with two incidence relations $h, t: V(D) \times
E(D) \rightarrow {0,1}$. We let $t(u,e) = 1$ iff $u$ is the tail of $e$ 
$e$ (i.e. the edge is directed from $u$ to another vertex $v$), and similarly 
$h(v,e) = 1$ iff $v$ is the head of the edge. Note this allows multiple 
edges from $u$ to $v$. We will assume that $|V(D)|$ and $|E(D)|$ are finite
throughout this paper, unless explicitly stated otherwise.
\item A {\it simple digraph} is a directed graph $G$ such that for 
all $u,v \in V(G)$, at most one edge from $u$ to $v$ appears in $E(G)$ 
(i.e. no parallel directed edges). 
\item In a digraph $G$, for $u,v \in V(G)$, 
the distance { $d(u,v)$} 
from $u$ to $v$ is the length of the shortest directed path from $u$ to $v$. 
We say $v$ is at {\it out-distance} $d(u,v)$ from $u$, and $u$ at 
{\it in-distance} $d(u,v)$ from $v$. 
\item For integers $j > 0$, { $N_j^+(v)$}
is the set of vertices at out-distance exactly $j$ from $v$, and {
$N_j^-(v)$} is the set of vertices at in-distance exactly $j$ from $v$. We may abbreviate
$N_1^+(v)$ and $N_1^-(v)$ to $N^+(v)$ and $N^-(v)$, respectively.
\item In a digraph $G$, { $\delta_G^+$ and $\delta_G^-$ }
denote the minimum out-degree and in-degree of $G$, respectively. For a given vertex $v$, { 
$\delta^+_G(v), \delta^-_G(v)$} denote the out-degree and in-degree of the vertex $v$. 
In cases where the graph being referenced is clear, we may write $\delta^+(v)$ and $\delta^-(v)$.
\end{itemize}

\section{The Caccetta-H\"{a}ggkvist Conjecture} 

\begin{conj} \label{CH} (L. Caccetta, R. H\"{a}ggkvist ~\cite{caccettahaggkvist}) 
Every simple $n$-vertex digraph with minimum out-degree at least 
$r$ has a cycle with length at most $\lceil \frac{n}{r} \rceil$.
\end{conj}

This can be restated in the following way: 
Let $A$ be an $n \times n$ $0$-$1$ matrix such that $a_{ij} = 1$ implies $a_{ji} \neq 1$ for all 
$i \neq j$. Let $a_{ii}=1$ for all $i$. If the sum of every row in $A$ is at least $r+1$, 
$A^{\lceil n/r \rceil}$ has trace greater than $n$.

\subsection{Partial Results} 
The C-H conjecture has been proved for: 
\begin{itemize}
\item $r = 2$ by Caccetta and H\"{a}ggkvist ~\cite{caccettahaggkvist}
\item $r = 3$ by Hamidoune ~\cite{hamidouneCH} 
\item $r = 4$ and $r = 5$ by Ho\'{a}ng and Reed ~\cite{hoangreed}
\item $r \leq \sqrt{n/2}$ by Shen ~\cite{shen:dm2000}. For the exact statement
of his result, see Theorem \ref{shensqrt}. This shows that for any given $r$, 
the number of counterexamples to the conjecture (if any) is finite.
\item Cayley graphs (which implies all vertex transitive graphs using 
coset representations) by Hamidoune ~\cite{hamidouneCayley}. 
This proof uses a lemma of Kemperman ~\cite{kemperman} (Lemma 
~\ref{kempermanlemma}). 
\end{itemize}

Also, Shen ~\cite{shen:ejc2003} proved that if $deg^+(u)+deg^+(v) \geq 4$ 
for all $(u,v) \in E(G)$, then $g \leq \lceil n/2 \rceil$, where $g$ denotes the girth of $G$. This is an average local outdegree 
version for the $r=2$ case of the Caccetta-H\"{a}ggkvist conjecture.

\subsection{Approximate Results I - Additive Constant}

Another approach is to show that if $\delta^+_G \geq r$, then there is 
a cycle of length at most $\frac{n}{r} + c$ for some small $c$. This has 
been proved for some values of $c$, as follows:

\begin{itemize}
\item $c = 2500$ by Chv\'{a}tal and Szemer\'{e}di ~\cite{chvatalszemeredi}
\item $c = 304$ by Nishimura ~\cite{nishimura}
\item $c = 73$ by Shen ~\cite{shen:gc2002}.
\end{itemize}

\subsection{Approximate Results II - Special Case n/3}

The case $r = n/2$ is trivial, but $r = n/3$ has received much attention. 
Research has sought the minimum constant $c$ such that $\delta^+_G \geq cn$ 
in an $n$-vertex simple digraph $G$ forces a 
directed cycle of length at most 3. The conjecture is that $c = 1/3$, and
the current results are:

\begin{itemize}
\item $c \leq (3-\sqrt{5})/2 = 0.382$ by Caccetta and H\"{a}ggkvist 
~\cite{caccettahaggkvist}
\item $c \leq (2\sqrt{6}-3)/5 = 0.3797$ by Bondy in a neat subgraph 
counting argument ~\cite{bondy}
\item $c \leq 3 - \sqrt{7} = 0.3542$ by ~\cite{shen:jct1998}
\end{itemize}

Similarly, Seymour, Graaf, and Schrijver ~\cite{seymourgraafschrijver}
asked for the minimum value of 
$\beta$ so that when the minimum in- and out-degrees of $G$ are at least 
$\beta n$, $G$ has a directed cycle of length at most $3$. They proved that 
$\beta \leq 0.3487$ and gave a formula relating $\beta$ and $c$. Shen applied 
this formula to his 1998 result to get a slight improvement to $\beta \leq
0.3477$ ~\cite{shen:jct1998}.

\section{Seymour's Second Neighborhood Conjecture}

This conjecture implies the special case of 
Caccetta-H\"{a}ggkvist when both in- and out-degrees are at least $n/3$,
and has received much attention of its own. 

\begin{conj} (Seymour) Any simple digraph with no loops or digons 
has a vertex $v$ whose second neighborhood is at least as big as its first neighborhood, 
i.e. $|N^+_2(v)| \geq |N^+(v)|$. 

\end{conj}

The following is known for Seymour's second neighborhood conjecture: 
\begin{itemize}
\item When G is a tournament, this is Dean's conjecture, and was proved by 
Fisher ~\cite{fisher} using probabilistic methods. There is also a
combinatorial proof by Havet and Thomass\'{e} ~\cite{havetthomasse}.
\item It was proved for digraphs with minimum outdegree $\leq 6$ by 
Kaneko and Locke ~\cite{kanekolocke}.
\item There is a vertex $v$ where $|N^+_2(v)| \geq \gamma|N^+(v)|$ and 
$\gamma = 0.657298...$ is the unique real root of $2x^3 + x^2 -1=0$. (Note
the conjecture is that $\gamma = 1$) by Chen, Shen, and Yuster 
~\cite{chenshenyuster}. 
They also claim a slight improvement to $\gamma = 0.67815$ (proof unpublished).
\item Godbole, Cole, and Wright ~\cite{godbolecolewright}
 showed that the conjecture holds for almost all digraphs.
\end{itemize}

\section{r-Regular Digraphs}

A digraph $G$ is {\it  $r$-regular } if every vertex $v$ has 
$\delta^+_G(v)=\delta^-_G(v) = r$.

\begin{conj}(Behzad, Chartrand, Wall ~\cite{bhc}) 
The minimum number of vertices in an $r$-regular digraph $G$ with 
girth $g$ is $r(g -1)+1$. 
\end{conj}

Behzad, Chartrand and Wall give an example achieving this by placing
 $r(g-1)+1$ vertices on a circle with each vertex having edges 
to the next $r$ vertices in clockwise order. The Caccetta-H\"{a}ggkvist 
Conjecture is a generalization of this earlier conjecture.

The Behzad-Chartrand-Wall conjecture was proved for the following special cases:
\begin{itemize}
\item $r = 2$ by Behzad ~\cite{behzad}
\item $r = 3$ by Bermond ~\cite{bermond}
\item Vertex-transitive graphs by Hamidoune ~\cite{hamidouneVT}
\item If $\delta^+_G \geq r$, then 
$g \leq 3\lceil \frac{n}{r}\ln(\frac{2+\sqrt{7}}{3})\rceil 
\approx \frac{1.312n}{r}$ by Shen ~\cite{shen:gc2002}.
\end{itemize}

\section{Related Results}

\begin{thm}\label{shensqrt} (Shen ~\cite{shen:dm2000}) For a digraph $G$ on $n$ vertices, 
if $\delta^+_G \geq r$, and $n \geq 2r^2 - 3r + 1$, then $G$ has a cycle of 
length at most $\lceil \frac{n}{r} \rceil$. 
\end{thm}

 In a graph $G$, for $u,v \in V(G)$, 
$\kappa(u,v)$ denotes
the maximum number of internally disjoint paths between $u$ and $v$. If $G$ is a digraph, 
$\kappa$ counts the maximum number of 
internally disjoint directed paths from $u$ to $v$.

In a graph $G$, for $u,v \in V(G)$, { $\lambda(u,v)$} 
denotes the maximum number of edge-disjoint paths between $u$ and $v$. 
If $G$ is a digraph, $\lambda$ counts the maximum number of edge-disjoint 
directed paths from $u$ to $v$.

\begin{thm}\label{Thom85} (Thomassen ~\cite{thomassen}) For all positive integers 
$r$, there is a digraph $D$ without digons with 
$\delta^+_D \geq r$ and $\delta^-_D \geq r$ such that: 
\begin{enumerate}
\renewcommand{\labelenumi}{\roman{enumi}.}
\item no vertex $v \in V(D)$ is contained in three openly disjoint circuits (that is, 
three circuits which pairwise share only $v$)
\item no edge $(x,y) \in E(D)$ has $\kappa(y,x) \geq 3$.
\end{enumerate}
\renewcommand{\labelenumi}{\arabic{enumi}.}
\end{thm}

\begin{thm} (Mader ~\cite{mader95}) For every integer $k \geq 0$, If $G$ has
$|V(G)| > k^2(k+1)$ and at most $k^2(k+1)$ vertices of $G$ 
have out-degree at most $k^3(k+1)$, then there are vertices $x \neq y$ 
so that $\kappa(x,y) > k$. 
\end{thm}

\subsection{Undirected Graph Theorems}

\begin{thm} (Mader ~\cite{mader71}) Every graph $G$ with $\delta_G \geq r$ contains 
vertices $x,y$ with $\kappa(x,y) \geq r$ when $r \geq 1$. 
\end{thm}

\begin{thm} (Mader ~\cite{mader74}) Every graph $G$ with $\delta_G \geq r$ contains 
$r+1$ vertices $v_1, \dots, v_{r+1}$ with $\lambda(v_i, v_j) \geq r$ for all $i \neq j$.
\end{thm}

\subsection{Additive Number Theory Results}

 Given an additive group $\Gamma$, and sets $A,B \subseteq \Gamma$, 
let $A+B$  
$ := \{a + b \;|\; a \in A, b \in B\}$, and  $A \hat{+} B$ 
$ := \{a + b \;|\; a \in A, b \in B, a\neq b\}$. Finally, for a positive integer $r$, let 
 $rB$  $:= \{b_1 + \cdots + b_h |$ all $ b_i \in B$, not 
necessarily distinct$\}$.

\begin{thm} (Cauchy ~\cite{cauchy} and Davenport 
~\cite{davenport},~\cite{davenport2}) Let $p$ be a prime, and 
$A,B \subseteq \Z/p\Z$ be nonempty. Then $|A+B| \geq \min(p, |A|+|B|-1)$.
\end{thm}

\begin{thm} (I. Chowla ~\cite{chowla}) Let $m$ be a positive integer, and 
$A,B \subseteq \Z/m\Z$ such that $0 \in B$ and $gcd(b,m) = 1$ for all nonzero $b \in B$.
Then $|A+B| \geq \min(m, |A|+|B|-1)$.
\end{thm}

\begin{thm} (Dias de Silva and Hamidoune ~\cite{desilvahamidoune}) 
\label{erdosheilbronn} The Erd\"{o}s-Heilbronn Conjecture:  
Let $A \subseteq \Z/p\Z$, with $p$ prime. 
Then $|A \hat{+} A| \geq \min(2|A|-3, p)$. 
\end{thm}

 For a multiplicative group $\Gamma$ and sets $A,B \subseteq \Gamma$, 
let { $AB$} $:= \{ab \;|\; a\in A, b\in B\}$.

\begin{lemma}\label{kempermanlemma} (Kemperman ~\cite{kemperman}) 
Given a group $\Gamma$ and finite 
non-empty subsets $A,B \subseteq \Gamma$, if $1 \in A, B$ but 
$(1,1)$ is the only 
pair $(a,b)$ with $a \in A, b \in B$ such that $ab=1$, then $|AB| \geq |A|
+ |B| -1$.  
\end{lemma}

 We say $G$ is a {\it layered digraph} if $G$ is a digraph with 
$V(G) = \cup_{i=0}^{h} V_i$ with $V_i \neq 
\emptyset$, and $V_i \cap V_j = \emptyset$ for all $i \neq j $, and $(u,v) \in E(G)$ implies
$u \in V_{i-1}, v \in V_i$ for some $i \in \{1,\dots,h\}$.

 A {\it Pl\"unnecke graph} is a layered digraph $G$ with the following 
two properties: 
\begin{enumerate}
\item If $u,v,w_1, \dots, w_k$ are vertices of $G$ with 
$(u,v), (v,w_1), \dots, (v,w_k) \in E(G)$, then there are distinct vertices 
$v_1, \dots, v_k$ so that $(u,v_i), (v_i, w_i) \in E(G)$ for $i=1,\dots,k$. 
\item If $v,w, u_1, \dots, u_k$ are vertices of $G$ with 
$(v,w), (u_1,v), \dots, (u_k,v) \in E(G)$, then there are distinct vertices 
$v_1, \dots, v_k$ so that $(u_i,v_i), (v_i, w) \in E(G)$ for $i=1,\dots,k$. 
\end{enumerate}

 Let $G$ be a digraph, and $X,Y$ nonempty subsets of $V(G)$. Then 
\\{\it Im(X,Y)} $:=\{y \in Y |$ there is a directed path from $X$ to $y\}$. 
The {\it magnification ratio D(X,Y)} is $$D(X,Y) := \min_{Z\subseteq X, Z \neq \emptyset}
\left\{\frac{|Im(Z,Y)|}{|Z|}\right\}.$$

\begin{thm} (Pl\"{u}nnecke ~\cite{plunnecke}) In a Pl\"{u}nnecke graph, let $D_i = D(V_0, V_i)$. Then 
$$D_1 \geq D_2^{\frac{1}{2}} \geq \cdots \geq D_h^{\frac{1}{h}}.$$
\end{thm}

The following are consequences of applying Pl\"{u}nnecke's Inequalities 
to a special graph created from subsets $A,B$ of a group:

\begin{thm} For sets $A,B \subset \Gamma$:
\begin{enumerate}
\item $|iB|^{\frac{1}{i}} \geq |hB|^{\frac{1}{h}}$ for all $0 \leq i < h$.  
\item If $|B| = k$, and $|B+B| \leq ck$, then $|hB| \leq c^hk$.
\item If $|A| = n$, and $|A+B| < cn$ then for all $k, \ell \in \Z^+$, we have that 
$|kB-\ell B| \leq c^{k+\ell}n$ where $kB-\ell B$ denotes the set of all elements 
expressable as $(b_1 + \cdots + b_k) - (b_1'+ \cdots + b_{\ell}')$ where all 
$b_i, b_i'$ are in $B$.
\end{enumerate}
\end{thm}

\begin{thm} (generalization of Erd\"{o}s-Heilbronn (Thm. \ref{erdosheilbronn})) 
Let $A,B \subseteq \Z/p\Z$, 
with $p$ prime and $|A| \neq |B|$. Let $C = A \hat{+} B$. 
Then $|C| \geq \min(|A| + |B| - 2, p)$.
\end{thm}

\section {Open Problems and Conjectures}

\subsection{Rainbow Conjectures}

\subsubsection{A Colored Generalization of Seymour's Second Neighborhood}

Given a digraph $G = (V,E)$ with each edge $e \in E$ having 
a set $S_e$ of labels in $\{1,2,\dots,k_G\}$, a {\it rainbow} 
structure $H$ in $G$ (such as a path or cycle) means one in which 
there is a way to assign each edge $e \in E(H)$ a label $\ell (e) \in S_e$ 
so that $ell (e) \neq \ell (f)$ for all edges $e \neq f$ in $E(H)$.

\begin{conj}\label{colorconj} (Seymour, Sullivan)
Let $G$ be a simple digraph on the vertex set $V$, and 
$E_1, \dots, E_k \subseteq E(G)$. 
Say an edge $e \in E(G)$ has label set $S_e \subseteq \{1,\dots,k\}$
where $i \in S_e$ if and only if $e \in E_i$. Finally, let $G_i = (V, E_i)$. 
\begin{enumerate}
\renewcommand{\theenumi}{\roman{enumi}}
\item There exists a rainbow (di)cycle in $G$ or
\item There exists a vertex $v$ such that $|\{w \,|\,$there exists a 
rainbow path from $v$ to $w\}| \geq  \sum_{i=1}^{k} \delta_{G_i}^+(v)$. 
\end{enumerate}
\end{conj}
\begin{proof}[Notes:] This is false if you require that the colors appear
in an increasing order (cyclic on the cycle). 
We have been able to show that this conjecture holds when
$G_1,\dots, G_k$ are Cayley graphs on a common group $\Gamma$ (using 
induction on Lemma ~\ref{kempermanlemma}), and when 
$\delta_{G_i}^+(v) \leq 1$ for all $v$ and
all $i$ except $i = 1$, where we allow the outdegrees to be unbounded 
(but finite). 
\end{proof}

\subsubsection{Implications of Conjecture \ref{colorconj}}
\begin{conj}
Seymour's Second Neighborhood Conjecture
\end{conj}
\begin{proof}[Notes:] 
To see this, for a digraph $H$, let $k=2, G=H$, and $E_1=E_2=E(H)$. 
\end{proof}

\begin{conj}
Caccetta-H\"{a}ggkvist Conjecture (general case)
\end{conj}
\begin{proof}[Notes:] 
For a digraph $H$, take $G=H$,
$k = \lceil \frac{n}{\delta^+_H} \rceil$, and $E_1=\dots=E_k=E(H)$. We must
get a rainbow cycle because the sum of the outdegrees at each vertex is
$\geq n$. This corresponds to a dicycle of length at most 
$\lceil \frac{n}{\delta^+_H} \rceil$ in $H$, as desired.
\end{proof}

\begin{conj}\label{secondinconj}
Any simple digraph with no loops or digons 
has a vertex $v$ such that $|N^+_2(v)| + |N^{+}(v)| \geq 2|N^{-}(v)|$.
\end{conj}
\begin{proof}[Notes:] 
Recall, for comparison,
SSN can be written as $|N^+_2(v)| + |N^{+}(v)| \geq 2|N^{+}(v)|$. 
To see how Conjecture
 \ref{colorconj} implies \ref{secondinconj}, 
take a simple digraph $H$ with no loops or digons, 
and set $G=H$ and $E_1 = E_2 = E(H)$. 
$G$ cannot have a rainbow cycle by definition of $H$. 
Define 
$N^{+*}_G(u) = \{$vertices you can reach by a rainbow path in $G$ from $u\}$, and 
$N^{-*}_G(u) = \{$vertices that have a rainbow path in $G$ to $u\}$. 
Let $E_3$ be the edges $\{(u,v)\;|\; v$ is not in the set $N^{-*}_G(u)\}$. 
Let $G' = H$, and have subsets $E_1, E_2,E_3 \subseteq E(G')$
giving rise to label sets $S'_e \subseteq \{1,2,3\}$ for $e \in G'$.
We can see from these definitions that for any $u$ in $V(G')$, 
$$\sum_{i=1}^{3} \delta_{G_i}^+(u) = 2|N^+(u)| + ((n-1) - |N^-(u)| - 
|N^-_2(u)|),$$
where all neighborhoods referenced on the RHS are in $H$. We can rewrite this as: 
$$\sum_{i=1}^{3} \delta_{G_i}^+(u) = ((n-1) - (|N^-(u)| + |N^-_2(u)|-2|N^+(u)|)).$$
Then $\sum_{i=1}^{3} \delta_{G_i}^+(u) \geq n$ whenever $|N^-(u)| + |N^-_2(u)| < 2|N^+(u)|.$ 
If this were true for all vertices $u$, then no vertex could have $|N^{+*}_{G'}(u)| \geq 
\sum_{i=1}^{3} \delta_{G_i}^+(u) \geq n$, so we must have a rainbow cycle in $G'$, by 
Conjecture \ref{colorconj}. However, by construction, since $G$ has no rainbow cycle, 
$G'$ has no rainbow cycle. Thus there is a vertex
$v \in V(H)$ so that $|N^-(u)| + |N^-_2(u)| \geq 2|N^+(u)|.$ If we reverse all edges in $H$, this
gives  $|N^+(u)| + |N^+_2(u)| \geq 2|N^-(u)|,$ as claimed. 
\end{proof}

\begin{conj}\label{forcecycleconj} (Seymour)
Under the hypotheses of Conjecture \ref{colorconj}, if $|V| = d$ and $\sum_{i=1}^{k} 
\delta_{G_i}(v) \geq d$ for all $v$, $G$ must have a rainbow cycle.
\end{conj}

\begin{proof}[Note:] This conjecture is false if $|V| =d$ is replaced by $|V|=d+1$. \end{proof}

\subsubsection{Other Conjectures Inspired by (or related to) Conjecture \ref{colorconj}}
If we believe Seymour's second neighborhood conjecture and Conjecture \ref{colorconj} 
(specifically \ref{secondinconj}), we might be led to ask if the following holds:
\begin{conj} (``Compromise Conjecture'')
Any simple digraph with no loops or digons 
has a vertex $v$ such that $|N^+_2(v)| \geq |N^{-}(v)|$.
\end{conj}

\begin{conj}
Any simple digraph with no loops or digons has
a vertex $v$ such that $|N^+_2(v)|+|N^+(v)| \geq 2\min(|N^{-}(v)|, |N^{+}(v)|)$.
\end{conj}

\begin{conj} Under the hypotheses of Conjecture \ref{colorconj}, if 
$\delta^+_{G_i} \geq r_i$, and 
$\sum_{i=1}^{t} r_i \geq |V|$, there is a rainbow cycle in $G$.
\end{conj}

\begin{conj} Under the hypotheses of Conjecture \ref{colorconj}, if 
$\sum_{i=1}^{t} \delta^+_{G_i}(v) \geq |V|$ for all vertices $v$, 
there is a rainbow cycle in $G$.
\end{conj}

\begin{conj} (Devos) Under the hypotheses of Conjecture ~\ref{colorconj}, if
there is no rainbow cycle in $G$ with strictly increasing edge labels, 
then the average number of vertices reachable from a fixed
vertex $v$ by label-increasing (possibly trivial) 
paths is at least $1+\sum_{i=1}^k \delta^+_{G_i}$.
\end{conj}
\begin{proof}[Note:] This can be proved when 
$G_i =$ Cayley$(\Gamma,A_i)$ for some group $\Gamma$, 
using Lemma ~\ref{kempermanlemma}
\end{proof}

\subsection{Second \& $K^{th}$ Neighborhood Conjectures}

\begin{conj} Is Seymour's Second Neighborhood true for 
locally finite digraphs? What if we just require the outdegrees to be finite?
\end{conj}

\begin{conj} (Thomass\'{e}) Let $G$ be a digraph with no directed cycle of length at most three. Then 
there is a vertex $v \in V(G)$ with $\delta^+(v)$ at most the number of non-neighbors of $v$.
\end{conj}

The following generalization of second neighborhood to $k$th neighborhood 
was taken from a (unpublished) paper of Serge Burckel: 
\begin{conj}
Any simple digraph with no directed cycles of length at most $k$
has a vertex $v$ such that $|N^+_k(v)| \geq |N^+_{k-1}(v)|$. 
\end{conj}

Serge Burckel also asked the following structural question: 
\begin{conj}
For any $k$ and any digraph $G$, define $G^*$ to be
those vertices with $|N^+_2(v)| \geq |N^+(v)|$. Then any vertex of 
out-degree $k$ is at distance at most $k$ from a vertex in $G^*$.
\end{conj}
\begin{proof}[Notes:]
He motivates this with the following remark: 
"If a vertex $x$ has one successor $y$, then if $y$ has no successors, 
$y \in G^*$, otherwise $x \in G^*$. This property seems to generalize for
any out-degree, and if it is true, is optimal considering 'pyramids' where
any vertex of out-degree $k$ is at distance exactly $k$ from the (unique)
solution."
His 'pyramids' are formed by placing one vertex, then two in the row beneath 
it, and so forth ($i$ vertices in row $i$) to form a triangle. The vertices
in row $i$ are then completely joined to all vertices in row $i-1$ for
$i \geq 2$. 
\end{proof}

\begin{conj} (Seymour \&/or Jackson) 
If $G$ is an Eulerian digraph with no loops or digons, then 
$$\sum_{v\in V(G)} |N^+_2(v)| \geq \sum_{v \in V(G)} |N^+(v)|.$$
\end{conj}

\begin{conj} (Thomass\'{e} \& Kral) Let $G$ be an Eulerian digraph on $n$ vertices, 
so that $|E(G)| \geq n^2/3$. Then $G$ has a directed cycle of length at most three.
\end{conj}

\subsection{Matrices}

\begin{conj} The following were all presented by Seymour, with no other attributions given:
For the following questions, matrices are assumed to be 
$n \times n$ $0$-$1$ matrices such that $a_{ij} = 1$ implies $a_{ji} \neq 1$, and
all diagonal elements equal to $1$.
\begin{enumerate}
\item Let $A_1, A_2, \dots, A_{\lceil n/r \rceil}$ be matrices (not necessarily distinct) so
that the row sums of $A_i$ are at least $r+1$ for all $i$. Does $A_1A_2\cdots A_{\lceil n/r \rceil}$ 
have trace $>n$? This is a special case of Conjecture \ref{forcecycleconj}.
\item Let $A_1, A_2, \dots, A_t$ be matrices (not necessarily distinct) so
$A_i$ has row sums at least $r_i + 1$ and 
$\sum_{i=1}^t r_i \geq n$. Does $A_1A_2\cdots A_t$ have trace $>n$? This is equivalent to Conjecture 
\ref{forcecycleconj}.
\item Form a digraph from the matrices $A_1, \dots, A_t$ by putting $t+1$ copies of the vertex set $V$ 
in a row, connecting copy $k$ of $v$ to copy $k+1$ of $v$ 
with a horizontal edge for all vertices $v$ and all $k$,
and then from copy $k$ to copy $k+1$ put in the edge from copy $k$ of $v_i$ to copy $k+1$ of $v_j$ 
precisely when $a_{ij} = 1$ in matrix $A_k$ for $k=1,\dots,t$. The question now is 
whether there exists a vertex $u$ so that there is a non-trivial
(i.e. not all horizontal edges) path from copy $1$ of $u$ to copy $t+1$ of $u$. 
\item If we ``squish'' all the bipartite graphs from the previous items so they live on a single copy of 
$V$, marking edges from copy $i$ to copy $i+1$ with color $i$, then we have the sum of the (colored) 
outdegrees at each vertex is at least $n$, and we're asking for a non-trivial rainbow cycle which has 
the colors appearing in increasing order.
\end{enumerate}
\end{conj}

 For a matrix $A$, the {\it spectral radius of A} is defined to be $\max \{ |\lambda| \; : 
\; \lambda$ an eigenvalue of $A\}$.

\begin{conj} (Charbit) Let $A$ be the adjacency matrix of a digraph $G$ with zero
on the diagonal and $a_{ij}=1$ if and only if the edge $(i,j) \in E(G)$. 
If the spectral radius of $A$ is at least $n/k$, then $G$ has a cycle of length 
$\leq k$.
\end{conj}

\subsection{Disjoint Cycles}

\begin{conj}\label{bermondthomassen} (Bermond-Thomassen) In a digraph $D$ with $\delta^+_D \geq 2k-1$, 
there are $k$ vertex disjoint cycles.
\end{conj}
\begin{proof}[Notes:] Open for $k \geq 3$, and the proof for $k=2$ by Thomassen is not intuitive.
\end{proof}

\begin{conj} (Ho\'{a}ng \& Reed) \label{hoangreed} If $G$ is a digraph with minimum outdegree
$r$, then there are directed cycles $C_1, \dots, C_r$ such that for all $\ell$,
$$|V(C_{\ell}) \cap (\cup_{i=1}^{\ell-1}V(C_i))| \leq 1.$$
\end{conj}

\begin{proof}[Notes:] The related conjecture that given minimum outdegree at 
least $r$, there should be a vertex $v$ with 
$r$ cycles through $v$ which are otherwise vertex-disjoint is false. 
A counterexample for $r=3$ was given by Thomassen in 1985 (see Theorem \ref{Thom85}
for his complete result). Adding the condition that the minimum 
indegree is also at least $r$ does not improve the veracity of the 
conjecture, though it is still open when indegree and outdegree are 
identically $r$ everywhere. 
\end{proof}

\subsection{Connectivity}

\begin{conj} (Hamidoune, 1981) Let $D$ be a digraph with $\delta^+_D \geq r$, 
$\delta^-_D \geq r$, and $r \geq 1$. Then there is an edge $(x,y)$ such that 
$\kappa(y,x) \geq r$.
\end{conj}

Note that a counterexample to the above conjecture 
for all $r \geq 3$ is given by Thomassen in Theorem~\ref{Thom85}.

\begin{conj} (Mader) If a digraph $D$ has $\delta^+_D \geq r$, then there are vertices
$x \neq y$ such that $\lambda(x,y) \geq r$. Also, there is an edge $(x,y)$ so that 
$\lambda(x,y) \geq r$.
\end{conj}

The first part of this conjecture was proven by Mader~\cite{mader85} 
for $\lambda(x,y) \geq r-1$.

\subsection{Weighted Versions}

\begin{conj}(Bollob\'{a}s \& Scott) \label{bollobasscott} Let $p:E(G) \rightarrow [0,1]$. If 
$\sum_{v\in N^+_G(u)} p(uv) \geq 1$ 
 and $\sum_{v\in N^-_G(u)} p(vu) \geq 1$ for all $u \in V(G)$, there is a directed cycle in $G$ of 
total weight $\geq 1$.
\end{conj}
\begin{proof}[Note:] There is a nice proof that there is a dipath of total weight at least $1$. 
\end{proof}

\begin{conj} (Zhang) Let $G$ be a digraph on $n$ vertices, 
and $f: E(G) \rightarrow \{0,1,\dots\}$. If $\sum_{e \in E^+(v)} f(e) \geq n/k$ 
for all $v \in V(G)$, then there is a directed
cycle $C$ such that $\sum_{e \in C} \frac{1}{f(e)} \leq k$.
\end{conj}
\begin{proof}[Counterexample:] (Charbit) Take the directed Cayley graph $G$ with group 
$G = \Z_8$, and generators $\{1,2\}$. Let the weight function $f$ be $1$ on 2-edges, and
$2$ on 1-edges (where r-edges are those coming from the generator r). Now let $k =3$. 
We can calculate $\sum_{e \in E^+(v)} f(e) = 4 \geq 8/3$ for all vertices $v$, but there
is no directed cycle with $\sum_{e \in C} \frac{1}{f(e)} \leq 3$.
\end{proof}

\subsection{Averaged Outdegree Conditions}

\begin{conj} If $D$ is a digraph on $n$ vertices with 
$$\sum_{v\in V(D)} \log(1+\frac{1}{\delta^+_D(v)}) \geq n\log(1+\frac{1}{r}),$$ 
then $D$ has a cycle of length at most $\lceil n/r \rceil$.
\end{conj}

\begin{proof}[Counterexample:] Take a transitive tournament of size $n-1$, and replace the edge from
the vertex of out-degree $n-2$ to the vertex of out-degree zero with a path of length 2 in the opposite
direction (thus increasing the number of vertices to $n$). 
\end{proof}

\begin{conj} (Shen ~\cite{shen:ejc2003}) 
Let $G$ be a digraph with $n$ vertices and minimum outdegree at least one.
If $\delta^+_G(v)+ \delta^+_G(u) \geq 2r$ for every edge $(u,v)$ in $G$, 
then the girth $g$ of $G$ is at most $\lceil n/r \rceil$.
\end{conj}

\begin{proof}[Note:] This was proved by Shen for $r=2$ in ~\cite{shen:ejc2003}.\end{proof}

\subsection{UnCategorized}

The following conjecture would imply that in a counterexample
for Caccetta-H\"{a}ggkvist for $r = n/3$, one can order the vertices so that at least 
$75\%$ of the edges go from left to right:

\begin{conj} (Chudnovsky, Seymour, Sullivan)
Let $G$ be a simple digraph with $k$ non-edges (unordered pairs $\{u,v\}$ 
where both $uv$ and $vu$ are not in $E(G)$). If $G$ has no directed cycle of 
length at most $3$, one can delete at most $k/2$ edges from $G$ and 
obtain a graph with no directed cycle. 
\end{conj}
\begin{proof}[Notes:]
We know that there are tight examples for this conjecture (transitive 
tournaments, $C_4$, and products of these). It is also known that a minimal 
counterexample has no source or sink vertex, and no directed cut. Kostochka 
recently proved all vertices in a minimal counterexample 
have at least $3$ and at most $(n-1)/2$ non-neighbors. He has an 
argument using these facts to show the conjecture for all $k \leq 14$. 
\end{proof}

\begin{conj} (Devos) 
For any digraph $G$ with no directed cycles of length at most three, 
there is a probability distribution $p$ on $V(G)$ such that 
at every vertex $v$, $p(N^+(v)) \geq p(N^-(v))$ and
$p(N^-_2(v)) \geq p(N^-(v))$, where $p(S) := \sum_{s\in S} p(s)$ for a set $S \subseteq V(G)$. 
\end{conj}
\begin{proof}[Notes:] For tournaments, such a distribution exists, and is unique. Its 
existence for a general digraph would imply Seymour's Second Neighborhood Conjecture, and actually 
it would suffice to just have a probability distribution $p$ so that 
$p(N^-_2(v)) \geq p(N^-(v))$ on average in $G$.
\end{proof}

 Given a digraph $G = (V,E)$, we say $F \subseteq E(G)$ is a 
{\it feedback arc set} if the digraph $G' = (V, E-F)$ has no directed cycles.

\begin{conj} (Lichiardopol's Conjecture) Every digraph $D$ has some minimal 
feedback arc set 
$F$ which contains a path of length $\delta^+_D$.
\end{conj}
\begin{proof}[Notes:] This implies Ho\'{a}ng \& Reed (\ref{hoangreed}), 
Caccetta-H\"{a}ggkvist (\ref{CH}), Bermond-Thomassen (\ref{bermondthomassen}), and 
Thomass\'{e} (\ref{thomassegirth}). 
\end{proof}

\begin{conj} (Mader, 1985) For all $k \in \Z^+$, there exists $r \in \Z^+$ such that
every $r$-out-regular digraph contains a subdivision of the transitive tournament 
on $k$ vertices.
\end{conj}
\begin{proof}[Notes:] This conjecture is known for $k = 3$ (where $r = 2$), and
$k=4$ ($r=3$ proven by Mader in 1996 ~\cite{mader96}). 
The existence of $r$ for $k =5$ is still not known,
though $r=6$ has been conjectured.
\end{proof}

\begin{conj} (Shen) For a digraph $G$ on $n$ vertices of girth $g$, define
$$t(G,r) = \sum_{u: \delta^+_G(u) < r} (r-\delta^+_G(u)).$$ If $\delta^+_G 
\geq 1$, then $n \geq r(g-1) + 1 - t(G,r)$. 
\end{conj}

\begin{conj} (Thomass\'{e}) If $G$ is a loopless digon-free digraph the maximum number
of induced directed 2-edge paths is $n^3/15 + \mathcal{O}(n^2)$. 

\end{conj}
\begin{proof}[Notes:] First, note that $n^3/15 +\mathcal{O}(n^2)$ can be obtained by 
substituting $C_4$'s inside $C_4$'s. Next, given a digraph $G$, let 
$\tau$ be the number of induced directed 2-paths, $\eta$ the number
of induced edges, and $\theta$ the number of cyclic triangles. Then 
if $|V(G)| = n$, 
$$\tau + \eta/3 + 2\theta = n^3/12 + \mathcal{O}(n^2) - 
\frac{1}{6}\sum_{v\in V(G)}\left( (\delta^+(v) - \delta^-(v))^2 + 
(\frac{n}{2}-\delta^+(v))^2 + (\frac{n}{2} -\delta^-(v))^2 \right).$$ 

F\"{u}redi rewrote the terms on the right hand side in terms of $\tau, \eta,$ and $\theta$, and showed
$$\tau \leq \frac{n^3}{12} + \mathcal{O}(n^2).$$ 
Bondy has a slight improvement of this result, proving that $$\tau \leq \frac{2n^3}{25}.$$
\end{proof}

Define an { \it $\alpha$/$\beta$-digraph } to be a digraph $D$ on 
vertex set $V$ and  edges defined by $\sigma_1,\dots,\sigma_{\beta}$ permutations 
(linear orders) on $V$ where the edge $(i,j) \in E \iff i<j$ in at least 
$\alpha$ of the $\sigma_i$. We say $G$ is a {\it majority digraph} if $\frac{\alpha}{\beta} > 2/3$.

\begin{conj} (Thomass\'{e} \& Charbit) Caccetta-H\"{a}ggkvist holds for 3/4-digraphs.
\end{conj}
\begin{proof}[Notes:] Majority digraphs have no cycles of length at most three. The class of 
3/4-digraphs is stable under substitution, and contains the circular interval graph on 
$3k+1$ vertices with outdegree $k$ going clockwise. The class includes all known 
extremal examples for C-H, yet this class of 3/4-digraphs seems manageably small. It is still open
whether or not 3/4 digraphs must have vertices $x$ of outdegree less than $n/3$. 
One can more generally ask if Caccetta-H\"{a}ggkvist holds for larger classes of majority digraphs.
\end{proof}

\begin{conj}(Thomass\'{e}) \label{thomassegirth} 
Every digraph $D$ has a path of length $\delta^+_D(g-1)$, where $g$ is the girth of $D$. 
\end{conj}
\begin{proof}[Note:] This is open for $g = 3$, and implies Caccetta-H\"{a}ggkvist.
\end{proof}

\begin{conj} (Thomass\'{e}) In a digraph $D$,  
$$\sum_{v\in V(D)} |\delta^+_D(v)-\delta^-_D(v)| + |\{(u,v) | d(u,v) \leq 2\}| \geq 2|E(D)| + 
|\{v| \delta^+_D(v) > \delta^-_D(v)\}|.$$
\end{conj}
\begin{proof}[Note:] This is exact for transitive tournaments.
\end{proof}

\begin{conj}(Thomass\'{e}) Let $G$ be a digraph on $n$ vertices with minimum outdegree at 
least $4n/15$ so that $G$ is maximal with no cycles of length at most three, and has no 
homogeneous set (in other words, G cannot be obtained by substitution). Then $G$ is a Cayley graph 
on $3k+1$ vertices with $1,\dots,k$ as generators for some $k$ (i.e. the circular interval graph 
with everyone joined to next $k$ clockwise).
\end{conj}
\begin{proof}[Note:] $\Z/15\Z$ with generators $S=\{1,2,4,8\}$ gives exactly $4n/15$.
\end{proof}

\subsection*{Acknowledgements}
This research was performed while on appointment as a U.S. Department of Homeland Security (DHS) Fellow under the DHS Scholarship and Fellowship Program, a program administered by the Oak Ridge Institute for Science and Education (ORISE) for DHS through an interagency agreement with the U.S Department of Energy (DOE). ORISE is managed by Oak Ridge Associated Universities under DOE contract number DE-AC05-06OR23100. All opinions expressed in this paper are the author's and do not necessarily reflect the policies and views of DHS, DOE, or ORISE.

\bibliographystyle{abbrv}
\bibliography{CaccettaHaggkvist}

\end{document}